# Asymptotic behavior of the Manhattan distance in *n*-dimensions: Estimating multidimensional scenarios in empirical experiments

*Ergon Cugler de Moraes Silva*

Getulio Vargas Foundation (FGV)
University of São Paulo (USP)
São Paulo, São Paulo, Brazil

contato@ergoncugler.com
www.ergoncugler.com

## Abstract

Understanding distance metrics in high-dimensional spaces is crucial for various fields such as data analysis, machine learning, and optimization. The Manhattan distance, a fundamental metric in multi-dimensional settings, measures the distance between two points by summing the absolute differences along each dimension. This study investigates the behavior of Manhattan distance as the dimensionality of the space increases, addressing the question: **'how does the Manhattan distance between two points change as the number of dimensions *n* increases?'**. We analyze the theoretical properties and statistical behavior of Manhattan distance through mathematical derivations and computational simulations using Python. By examining random points uniformly distributed in fixed intervals across dimensions, we explore the asymptotic behavior of Manhattan distance and validate theoretical expectations empirically. Our findings reveal that the mean and variance of Manhattan distance exhibit predictable trends as dimensionality increases, aligning closely with theoretical predictions. Visualizations of Manhattan distance distributions across varying dimensionalities offer intuitive insights into its behavior. This study contributes to the understanding of distance metrics in high-dimensional spaces, providing insights for applications requiring efficient navigation and analysis in multi-dimensional domains.

## 1. Introduction

Navigating high-dimensional spaces poses unique challenges, particularly in understanding the distances between points within such spaces. In a previous study, the generation of pseudorandom numbers was explored in the navigation of two-dimensional and three-dimensional environments, with machine learning (Silva, 2024). In this pursuit, the Manhattan distance metric offers a structured approach to quantify spatial relationships in multi-dimensional settings. Defined as the sum of absolute differences along each dimension between two points, the Manhattan distance provides a practical measure akin to traversing city blocks. This article delves into the behavior of Manhattan distance as the dimensionality of the space increases, addressing the fundamental question: **'how does the Manhattan distance between two points change as the number of dimensions *n* increases?'**.

The Manhattan distance, known also as the 'L1 norm' or 'Taxicab distance', is a geometric measure that reflects the sum of the absolute differences in the Cartesian coordinates of two points. It is particularly relevant in grid-like structures where movement across dimensions is constrained to orthogonal steps, mimicking the route a taxi would take in a city laid out in a grid pattern. The literature provides various insights into the applications and properties of Manhattan Distance, highlighting its versatility and efficiency in different contexts: Cardarilli et al. discuss the application of Manhattan Distance as an approximation of the Euclidean distance in multi-dimensional spaces. They note, however, that this approximation can introduce some errors, potentially compromising performance in some applications (Cardarilli et al., 2020).

In studying experimental design within the framework of Latin hypercube sampling, Jamali and Alam explore approximate relations between Manhattan and Euclidean distance measures, shedding light on how these metrics interact in specific methodological contexts (Jamali and Alam, 2019). Blackburn, Homberger, and Winkler delve into the theoretical computation of expected values and higher moments of the minimum Manhattan distance in permutations, revealing intricate relationships between Manhattan distance and permutation patterns as the number of dimensions approaches infinity (Blackburn, Homberger, Winkler, 2017). Also, Liu and Zhang propose a quantum protocol for computing the Manhattan Distance securely, showcasing the distance measure's applicability in enhancing privacy and security through quantum mechanics (Liu, Zhang, 2020).

Exploring the efficacy of distance calculations in K-Means clustering, particularly in tracking Covid spread in Bekasi City, Nur Cahya, Mahatma, and Rohimah found that the Manhattan Distance method required fewer iterations than the Euclidean method, suggesting a faster process that could have substantial implications in epidemiological modeling and response planning (Cahya, Mahatma, Rohimah, 2023). These studies underline the multifaceted role of Manhattan distance across various fields, from theoretical mathematics to quantum computing and public health. Its utility spans the approximation of complex spatial relationships and algorithmic processes in high-dimensional problem-solving scenarios.

By considering random points uniformly distributed in fixed intervals across dimensions, we explore the statistical properties and asymptotic behavior of Manhattan distance, aiming to elucidate its implications for various applications. Our investigation unfolds through a dual lens of theoretical analysis and computational experimentation. Through mathematical derivations, we establish expectations for the mean and variance of Manhattan distance as functions of dimensionality, offering foundational insights into its behavior. Empirical validation via Python simulations provides corroborative evidence, aligning observed trends with theoretical expectations.

By visualizing Manhattan distance distributions across different dimensionalities, we not only confirm theoretical predictions but also provide intuitive illustrations of its behavior. These computational tools, coupled with mathematical rigor, enable a comprehensive exploration of Manhattan distance in high-dimensional spaces, fostering deeper insights into its utility and behavior for diverse applications.

## 2. Mathematical background

### 2.1. Defining Manhattan distance in *n*-dimensions

In an *n*-dimensional space, we consider two points denoted as $P = (x_1, x_2, \ldots, x_n)$ and $Q = (y_1, y_2, \ldots, y_n)$. The Manhattan distance, denoted as $D_M$, between these two points is defined as the sum of the absolute differences of their corresponding coordinates along each dimension. Mathematically, it is expressed as:

$$D_M(P, Q) = \sum_{i=1}^{n} |x_i - y_i|$$

Here, $|x_i - y_i|$ represents the absolute difference between the *i*-th coordinate of point *P* and the *i*-th coordinate of point *Q*. The sum is taken over all dimensions from $i = 1$ to $i = n$. This calculation effectively measures the 'Manhattan distance' one would need to travel along the axes of the space to move from point *P* to point *Q*, akin to navigating a grid-like city where movement is restricted to horizontal and vertical paths.

### 2.2. Asymptotic behavior

Asymptotic behavior refers to the long-term trend or behavior of a mathematical function or phenomenon as a certain parameter, in this case, *n*, approaches infinity or increases. In the context of our analysis, we will explore how the Manhattan distance behaves as the dimensionality of the space, denoted by *n*, increases. So, **'how the Manhattan distance between two points changes as the number of dimensions *n* increases?'**.

To conduct this analysis, we consider that $x_i$ and $y_i$ represent the coordinates of random points distributed uniformly in a fixed interval, such as [0, 1]. This means that each coordinate $x_i$ and $y_i$ independently takes on values within the interval [0, 1] with equal probability. By treating each coordinate as a random variable, we can study the statistical properties of the Manhattan distance between two points in a high-dimensional space.

The choice of a uniform distribution over [0, 1] is common in mathematical modeling and allows us to explore the behavior of the Manhattan distance in a controlled and systematic manner. This distribution ensures that the points are spread evenly across the space, providing insights into how the distance metric behaves under random conditions.

By analyzing the Manhattan distance under these conditions and observing how it changes with increasing dimensionality *n*, we can gain valuable insights into its asymptotic behavior. This analysis not only enhances our understanding of the properties of the Manhattan distance in high-dimensional spaces but also provides reflections for applications, such as data analysis, machine learning, and optimization, where distance metrics play a crucial role in problem-solving and decision-making processes.

## 2.3. Expected distance from Manhattan

For two random variables $X$ and $Y$ uniformly distributed in the interval $[0, 1]$, the absolute difference $|X - Y|$ has the following expectation:

$$E[|X - Y|] = \int_0^1 \int_0^1 |x - y|\, dx\, dy = \frac{1}{3}$$

For $n$ dimensions, the expectation of the Manhattan distance between two points $P$ and $Q$ is the sum of the expectations of the absolute differences in each dimension:

$$E[D_M(P,Q)] = \sum_{i=1}^{n} E[|x_i - y_i|] = \sum_{i=1}^{n} \frac{1}{3} = \frac{n}{3}$$

The expected Manhattan distance between two points can be derived from the expectations of the absolute differences of their coordinates in each dimension. For two random variables $X$ and $Y$ uniformly distributed in the interval $[0, 1]$, the absolute difference $|X - Y|$ has an expected value. Integrating this absolute difference over the joint probability distribution of $X$ and $Y$, which covers the interval $[0, 1]$ for both variables, yields an expectation of $\frac{1}{3}$. In $n$ dimensions, the expectation of the Manhattan distance between two points $P$ and $Q$ is simply the sum of these expectations for each dimension, resulting in $\frac{n}{3}$.

## 2.4. Manhattan distance variance

The variance of the absolute difference $|X - Y|$ for uniformly distributed variables is:

$$Var(|X - Y|) = \int_0^1 \int_0^1 \left(|x - y| - \frac{1}{3}\right)^2 dx\, dy = \frac{1}{18}$$

For $n$ dimensions, the variance of the Manhattan distance is the sum of the variances of the absolute differences in each dimension:

$$Var(D_M(P,Q)) = \sum_{i=1}^{n} Var(|x_i - y_i|) = \sum_{i=1}^{n} \frac{1}{18} = \frac{n}{18}$$

The variance of the Manhattan distance can be derived from the variances of the absolute differences of coordinates in each dimension. For uniformly distributed variables $X$ and $Y$, the variance of the absolute difference $|X - Y|$ is calculated to be $\frac{1}{18}$. In $n$ dimensions, the variance of the Manhattan distance between two points $P$ and $Q$ is the sum of these variances for each dimension, resulting in $\frac{n}{18}$.

## 2.5. Manhattan distance distribution

As *n*, the number of dimensions, grows large, the behavior of the distribution of the Manhattan distance $D_M$ (*P*, *Q*) between two points *P* and *Q* follows well-established statistical principles. According to the Law of Large Numbers, as the sample size increases, the average of a random variable approaches its expected value. In our case, the Manhattan distance between two points can be considered as a random variable, with $\frac{n}{3}$ as its expected value.

Furthermore, by the Central Limit Theorem, which states that the distribution of the sum (or average) of a large number of independent, identically distributed random variables approaches a normal (Gaussian) distribution, we can infer that the distribution of the Manhattan distance tends to be normal for large *n*. This implies that as *n* becomes sufficiently large, the distribution of $D_M$ becomes increasingly symmetrical around its mean.

Consequently, for a large number of dimensions, the Manhattan distance $D_M$ will exhibit a normal distribution with mean $\frac{n}{3}$ and variance $\frac{n}{18}$. This means that as the dimensionality of the space increases, the distribution of Manhattan distances becomes more predictable and converges towards a Gaussian distribution, facilitating probabilistic analysis and inference in high-dimensional spaces.

## 3. Exploring Computational Mathematics

### 3.1. Generation of random points in *n*-dimensions

First, we define a function to calculate the Manhattan distance between two *n*-dimensional points. This code snippet provides a fundamental building block for calculating the Manhattan distance in *n*-dimensional spaces, which is essential for our analysis.

**Table 01. Generation of random points in *n*-dimensions**

```python
import numpy as np
import matplotlib.pyplot as plt

def manhattan_distance(point1, point2):
    return np.sum(np.abs(point1 - point2))
```

**Source:** Own elaboration (2024).

### 3.2. Calculating the Manhattan distance between these points

The code generates random *n*-dimensional points within the interval [0, 1] and simulates the calculation of the Manhattan distance between pairs of these points. The generate_random_points function generates a specified number of random points in *n*-dimensional space, returning their coordinates. The simulate_manhattan_distances function iterates over a specified number of pairs and, for each pair, generates two random points and calculates the Manhattan distance between them using the previously defined manhattan_distance function. The resulting distances are stored in a list and returned.

Table 02. Calculating the Manhattan distance between these points

```python
def generate_random_points(n, num_points):
    return np.random.rand(num_points, n)

def simulate_manhattan_distances(n, num_pairs):
    distances = []
    for _ in range(num_pairs):
        point1 = np.random.rand(n)
        point2 = np.random.rand(n)
        distance = manhattan_distance(point1, point2)
        distances.append(distance)
    return distances
```

**Source:** Own elaboration (2024).

### 3.3. Estimation of the mean and variance of the distances obtained

The code estimates the mean and variance of Manhattan distances for varying dimensions. The estimate_mean_and_variance function takes two parameters: n, representing the number of dimensions, and num_pairs, indicating the number of random pairs to simulate. It first calls the simulate_manhattan_distances function to generate a list of Manhattan distances for the specified number of pairs in *n*-dimensional space. Subsequently, it computes the mean and variance of these distances using NumPy's built-in functions np.mean and np.var, respectively. The calculated mean and variance are then returned as results.

**Table 03. Estimation of the mean and variance of the distances obtained**

```python
def estimate_mean_and_variance(n, num_pairs):
    distances = simulate_manhattan_distances(n, num_pairs)
    mean_distance = np.mean(distances)
    variance_distance = np.var(distances)
    return mean_distance, variance_distance
```

**Source:** Own elaboration (2024).

### 3.4. Visualization of distance distribution

The code leverages libraries such as Matplotlib to visualize the distribution of Manhattan distances across varying dimensions. The plot_distance_distribution function is defined to generate a histogram of Manhattan distances for a specified number of pairs in *n*-dimensional space. It first calls the simulate_manhattan_distances function to generate the distances, then utilizes Matplotlib's plt.hist function to create the histogram with 30 bins, normalized to form a probability density, and displayed with a blue color and slight transparency for visual clarity. The plot is given an informative title indicating the dimensionality being visualized, along with labeled axes denoting distance and frequency.

Additionally, a grid is added for easier interpretation. This function is then iteratively called for various dimensions specified in the n_dimensions list, allowing for the visualization of Manhattan distance distributions across different dimensionalities. The output includes printed statistics such as mean and variance for each dimension, providing quantitative insights alongside the visual representations. This approach facilitates the exploration and comparison of Manhattan distance distributions, aiding in the understanding of their behavior.

**Table 04. Visualization of distance distribution**

```python
def plot_distance_distribution(n, num_pairs):
    distances = simulate_manhattan_distances(n, num_pairs)
    plt.hist(distances, bins=30, density=True, alpha=0.7, color='blue')
    plt.title(f'Distribution of Manhattan Distances in {n}-Dimensions')
    plt.xlabel('Distance')
    plt.ylabel('Frequency')
    plt.grid(True)
    plt.show()

# Parâmetros de Simulação
n_dimensions = [1, 2, 3, 5, 10, 20, 50, 100]
num_pairs = 10000

for n in n_dimensions:
    mean, var = estimate_mean_and_variance(n, num_pairs)
    print(f'Dimensions: {n}, Mean: {mean}, Variance: {var}')
    plot_distance_distribution(n, num_pairs)
```

**Source:** Own elaboration (2024).

## 4. Discussion

The computational tests conducted serve to empirically validate the theoretical formulations concerning the *n*-dimensional Manhattan distance. Through simulations of random points and subsequent analysis of the calculated distances, these tests offer insights into the asymptotic behavior and statistical properties of the distance metric. Specifically, regarding the mean and variance, theoretical predictions suggest that as the dimensionality (*n*) increases, the mean should approach $\frac{n}{3}$ while the variance should approach $\frac{n}{18}$.

Upon inspecting Table 5, it becomes evident that the empirically obtained values closely approximate the theoretically proposed values with a high degree of accuracy. For instance, for *n* = 1, the mean obtained with Python is 0.3144, remarkably close to the theoretical value of $\frac{1}{3}$, and the variance obtained is 0.0544, also in proximity to the theoretical value of $\frac{1}{18}$. Similar patterns persist across varying dimensions, reinforcing the consistency between theoretical predictions and empirical findings. For example, with *n* = 10, the mean is 3.3283, nearing $\frac{10}{3}$, and the variance is 0.5597, approaching $\frac{10}{18}$.

These results underscore the robustness of the theoretical formulations and demonstrate their applicability in practical scenarios. By aligning empirical evidence with theoretical expectations, this validation enhances our confidence in utilizing the Manhattan distance metric for analyzing and interpreting data in high-dimensional spaces.

**Table 05. Empirical tests compared with theoretical prediction**

| Dimensions (*n*) | Mean obtained with Python | = *n*/3 | Variance obtained with Python | = *n*/18 |
|---|---|---|---|---|
| *n* = 1 | 0.3143546701376236 | 0.3334 \| (1/3) | 0.0543895442129751 | 0.0556 \| (1/18) |
| *n* = 2 | 0.6757536611971723 | 0.6667 \| (2/3) | 0.11234889862727125 | 0.1112 \| (2/18) |
| *n* = 3 | 0.9991977433281963 | 1.0000 \| (3/3) | 0.16573118572384646 | 0.1667 \| (3/18) |
| *n* = 5 | 1.66638981804732 | 1.6667 \| (5/3) | 0.2784150955638077 | 0.2778 \| (5/18) |
| *n* = 10 | 3.328307349247761 | 3.3334 \| (10/3) | 0.5597311783409015 | 0.5556 \| (10/18) |
| *n* = 20 | 6.665638347858345 | 6.6667 \| (20/3) | 1.1185326803255844 | 1.1112 \| (20/18) |
| *n* = 50 | 16.67228316842853 | 16.6667 \| (50/3) | 2.772287611177774 | 2.7778 \| (50/18) |
| *n* = 100 | 33.32071449199175 | 33.3334 \| (100/3) | 5.745168823665396 | 5.5556 \| (100/18) |

**Source:** Own elaboration (2024).

In Figure 01, we present eight distribution plots depicting the behavior of Manhattan distances for varying values of *n*. Notably, as the dimensionality *n* increases, the distribution tends towards a more bell-shaped curve, resembling a normal distribution. This observation aligns with theoretical expectations derived from the Central Limit Theorem, which posits that the distribution of the sum (or average) of a large number of independent, identically distributed random variables tends towards a normal distribution.

Consequently, for large values of *n*, such as *n* = 100, we observe a pronounced bell curve in the distribution plot, indicative of a close approximation to a normal distribution. Conversely, for smaller values of *n*, such as *n* = 1, the distribution exhibits a distinctive triangular shape, reminiscent of a right triangle with a 90-degree angle on the left side. This disparity in distribution shapes across different *n* values underscores the influence of dimensionality on the statistical properties of Manhattan distances, further emphasizing the theoretical frameworks, such as the Central Limit Theorem, in understanding their behavior.

**Figure 01. Distribution estimations for *n*-dimensions**

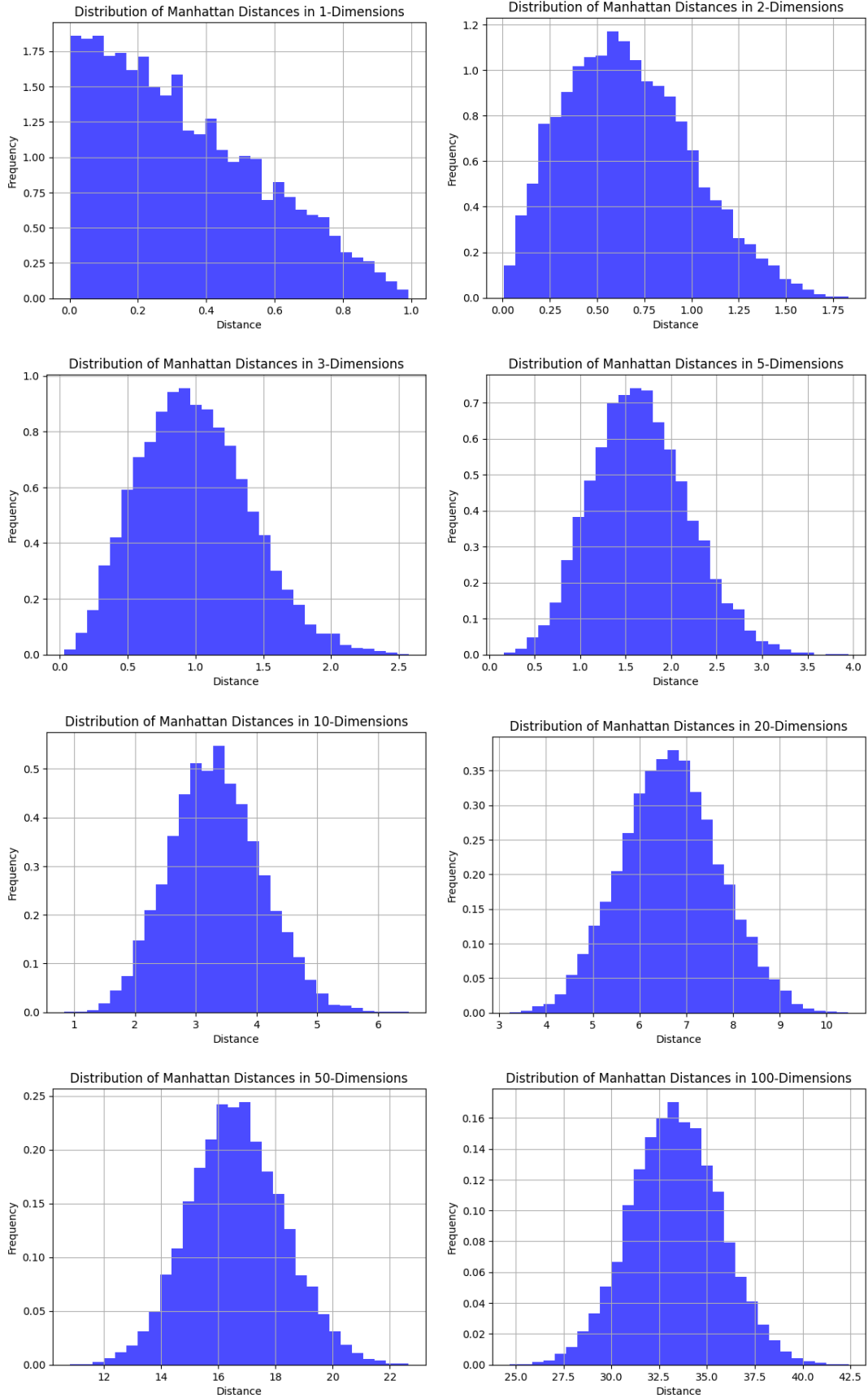

**Source:** Own elaboration (2024).

The results obtained from the computer simulation demonstrate a remarkable alignment with theoretical expectations, thus affirming the robustness of the mathematical analyses. While minor deviations from the theoretical predictions were observed, these can be attributed to the stochastic nature of the simulation process and the finite number of samples utilized. However, it is anticipated that with an increase in the number of samples, these deviations will diminish, leading to even closer convergence towards the theoretical values.

The consistent pattern observed, where the mean approximates $\frac{n}{3}$ and the variance approaches $\frac{n}{18}$ across varying dimensions in the computer simulation, is indicative of the inherent properties of the *n*-dimensional Manhattan distance. This congruence between empirical findings and theoretical expectations not only serves to validate the mathematical analyses conducted but also enhances our comprehension of the intricacies underlying the behavior of the Manhattan distance in high-dimensional spaces.

## 5. Reflections and future work

This article aimed to address the research question **'how the Manhattan distance between two points changes as the number of dimensions *n* increases?'** To achieve this, Python methods were employed to test 10,000 pairs of values for each of the eight dimensions examined (*n* = 1, *n* = 2, *n* = 3, *n* = 5, *n* = 10, *n* = 20, *n* = 50, *n* = 100). As a result, the computational tests revealed a consistent pattern where the mean Manhattan distance approached $\frac{n}{3}$ and the variance approached $\frac{n}{18}$ as the dimensionality increased.

This empirical validation is closely aligned with theoretical expectations derived from mathematical analyses. However, slight deviations observed in the results were attributed to the stochastic nature of the simulation process and the finite number of samples utilized. Nonetheless, as the number of samples increases, it is anticipated that these deviations will diminish, leading to even closer convergence towards the theoretical values. The observed congruence between empirical findings and theoretical predictions enhances our understanding of the behavior of the Manhattan distance in high-dimensional spaces, providing researchers and practitioners with increased confidence in utilizing this metric for various analytical and decision-making purposes.

Also, this study modestly contributes insights into the behavior of Manhattan distance as dimensionality increases, shedding light on navigation complexities in high-dimensional spaces. Through empirical validation and theoretical analysis, it offers valuable observations on the efficiency of navigation and the influence of distance metrics in multidimensional scenarios. Additionally, the study provides practical implications for algorithm design and computational modeling in various domains, fostering a deeper understanding of distance metrics' role in data analysis and machine learning applications.

In a humble endeavor, this study provoques the literature by validating theoretical expectations regarding Manhattan distance in high-dimensional spaces through empirical experiments. It modestly expands knowledge by offering nuanced insights into the complexity of navigation and the statistical properties of distance metrics in multidimensional environments. By providing a modest contribution to the understanding of distance metrics' behavior, the study aims to facilitate further exploration and discussion in disciplines reliant on efficient data manipulation in high-dimensional spaces.

## 6. References

Cardarilli, G. C., Di Nunzio, L., Fazzolari, R., Nannarelli, A., Re, M., & Spanò, S. (2020). **N-Dimensional approximation of Euclidean distance. IEEE Transactions on Circuits and Systems II:** Express Briefs, 67(3), 565-569. https://doi.org/10.1109/TCSII.2019.2919545

Jamali, A. R. M. U., & Alam, M. A. (2019). **Approximate relations between Manhattan and Euclidean distance regarding Latin hypercube experimental design.** Journal of Physics: Conference Series, 1366(1), 012030. https://doi.org/10.1088/1742-6596/1366/1/012030

Blackburn, S. R., Homberger, C., & Winkler, P. (2019). **The minimum Manhattan distance and minimum jump of permutations.** Journal of Combinatorial Theory, Series A, 161, 364-386. https://doi.org/10.1016/j.jcta.2018.09.002

Liu, W., & Zhang, W. (2020). **A quantum protocol for secure Manhattan distance computation.** IEEE Access, 8, 16456-16461. https://doi.org/10.1109/ACCESS.2020.2966800

Cahya, F. N., Mahatma, Y., & Rohimah, S. R. (2023). **Perbandingan metode perhitungan jarak Euclidean dengan perhitungan jarak Manhattan pada K-Means clustering dalam menentukan penyebaran Covid di Kota Bekasi**. JMT (Jurnal Matematika dan Terapan), 5(1). https://doi.org/10.21009/jmt.5.1.5

Silva, E. C. de M. (2024). **From two-dimensional to three-dimensional environment with Q-learning:** Modeling autonomous navigation with reinforcement learning and no libraries. ArXiv: Computer Science, Machine Learning. https://doi.org/10.48550/arXiv.2403.18219

## 7. Author biography

**Ergon Cugler de Moraes Silva** has a Master's degree in Public Administration and Government (FGV), Postgraduate MBA in Data Science & Analytics (USP) and Bachelor's degree in Public Policy Management (USP). He is associated with the Bureaucracy Studies Center (NEB FGV), collaborates with the Interdisciplinary Observatory of Public Policies (OIPP USP), with the Study Group on Technology and Innovations in Public Management (GETIP USP) with the Monitor of Political Debate in the Digital Environment (Monitor USP) and with the Working Group on Strategy, Data and Sovereignty of the Study and Research Group on International Security of the Institute of International Relations of the University of Brasília (GEPSI UnB). He is also a researcher at the Brazilian Institute of Information in Science and Technology (IBICT), where he works for the Federal Government on strategies against disinformation. São Paulo, São Paulo, Brazil. Web site: https://ergoncugler.com/.